\documentclass[a4paper,12pt]{amsart}

\usepackage{amsmath}
\usepackage{amsthm}
\usepackage{amsxtra}
\usepackage{amssymb}

\DeclareMathAlphabet{\mathds}{U}{dsrom}{m}{n}
\DeclareMathAlphabet{\mathsc}{U}{rsfs}{m}{n}

\theoremstyle{definition}
\newtheorem{DSBL-dfn}{Definition}
\theoremstyle{plain}
\newtheorem{DSBL-lem}[DSBL-dfn]{Lemma}
\newtheorem{DSBL-prp}[DSBL-dfn]{Proposition}
\newtheorem{DSBL-thm}[DSBL-dfn]{Theorem}

\begin{document}

\title{The Differential Structure of the Brieskorn Lattice}
\author{Mathias Schulze}
\address{M. Schulze, Department of Mathematics, D-67653 Kaiserslautern}
\email{mschulze@mathematik.uni-kl.de}

\begin{abstract}
The Brieskorn lattice $H''$ of an isolated hypersurface singularity with Milnor number $\mu$ is a free $\mathds{C}\{\!\{s\}\!\}$-module of rank $\mu$ with a differential operator $t=s^2\partial_s$.
Based on the mixed Hodge structure on the cohomology of the Milnor fibre, M. Saito constructed $\mathds{C}\{\!\{s\}\!\}$-bases of $H''$ for which the matrix of $t$ has the form $A=A_0+A_1s$.
We describe an algorithm to compute the matrices $A_0$ and $A_1$.
They determine the differential structure of the Brieskorn lattice, the spectral pairs and Hodge numbers, and the complex monodromy of the singularity.
\end{abstract}

\maketitle

\section{The Milnor Fibration}

Let $f:(\mathds{C}^{n+1},\underline{0})\longrightarrow(\mathds{C},0)$ be a holomorphic function germ with an isolated critical point and Milnor number $\mu=\dim_\mathds{C}\mathds{C}\{\underline{x}\}/\langle\underline{\partial}(f)\rangle$ where $\underline{x}=x_0,\dots,x_n$ is a complex coordinate system of $(\mathds{C}^{n+1},0)$ and $\underline{\partial}=\partial_{x_0},\dots,\partial_{x_n}$.
By the finite determinacy theorem, we may assume that $f\in\mathds{C}[\underline{x}]$.
By E.J.N~Looijenga \cite[2.B]{DSBL-Loo84}, for a good representative $f:X\longrightarrow T$ where $T\subset\mathds{C}$ is an open disk at the origin, the restriction $f:X'\longrightarrow T'$ to $T'=T\backslash\{0\}$ and $X'=X\backslash f^{-1}(0)$ is a $\mathsc{C}^\infty$ fibre bundle unique up to diffeomorphism, the Milnor fibration.
By J.~Milnor \cite[6.5]{DSBL-Mil68}, the general fibre $X_t=f^{-1}(t)$, $t\in T'$, is homotopy equivalent to a bouquet of $\mu$ $n$-spheres and, in particular, its reduced cohomology is $\widetilde{\mathrm{H}}^k(X_t)\cong\delta_{k,n}\mathds{Z}^\mu$ where $\delta$ is the Kronecker symbol.
Since $T'$ is locally contractible, the $n$-th cohomologies $\mathrm{H}(U)=\mathrm{H}^n(X_U)$ of $X_U=f^{-1}(U)$ form a locally free $\mathds{Z}$-sheaf of rank $\mu$ and $\mathrm{H}_\mathds{C}=\mathrm{H}\otimes_\mathds{Z}\mathds{C}$ is a complex local system of dimension $\mu$.
Hence, the sheaf of holomorphic sections $\mathsc{H}=\mathrm{H}\otimes_\mathds{Z}\mathsc{O}_{T'}$ of $\mathrm{H}_\mathds{C}$ is a locally free $\mathsc{O}_{T'}$-sheaf of rank $\mu$, the cohomology bundle.
By P.~Deligne \cite[2.23]{DSBL-Del70}, there is a natural flat connection $\nabla:\mathsc{H}\longrightarrow\mathsc{H}\otimes_{\mathsc{O}_{T'}}\Omega_{T'}^1$ on $\mathsc{H}$ with sheaf of flat sections $\mathrm{H}=\ker(\nabla)$, the Gauss-Manin connection.

\section{The Monodromy Representation}

Let $t$ be a complex coordinate of $T\subset\mathds{C}$, $i:T'\longrightarrow T$ the canonical inclusion, and $u:T^\infty\longrightarrow T'$ the universal covering of $T'$ defined by $u(\tau)=\exp(2\pi\mathrm{i}\tau)$ for a complex coordinate $\tau$ of $T^\infty\subset\mathds{C}$.
Then the covariant derivative $\nabla_{\partial_t}$ of $\nabla$ along $\partial_t$ induces a differential operator $\partial_t$ on $i_*\mathsc{H}$ and the pullback $f^\infty:X^\infty=X'\times_{T'}T^\infty\longrightarrow T^\infty$ is a $\mathsc{C}^\infty$ fibre bundle with $X^\infty_\tau=X_{u(\tau)}$, the (canonical) Milnor fibre.
Since $T^\infty$ is contractible, the $n$-th cohomologies $H(U)=\mathrm{H}^n(X^\infty_U)$ of $X^\infty_U=(f^\infty)^{-1}(U)$ form a free $\mathds{Z}$-sheaf of rank $\mu$ and $u_*H$ is the sheaf of multivalued sections of $\mathrm{H}$.
Lifting closed paths in $T'$ along sections of $\mathrm{H}$ defines the monodromy representation $\pi_1(T',t)\longrightarrow\mathrm{Aut}(\mathrm{H}_t)$ on $\mathrm{H}_t$ inducing the monodromy representation $\pi_1(T')\longrightarrow\mathrm{Aut}(H)$ on the cohomology $H$ of the Milnor fibre.
The image $\mathrm{M}$ of the counterclockwise generator of $\pi_1(T')$ is called the monodromy operator and fulfills $\mathrm{M}(s)(\tau)=s(\tau+1)$ for $s\in H$.
The sheaf $\mathrm{H}$ is determined by the monodromy representation up to isomorphism.
The following well known theorem is due to E.~Brieskorn \cite[0.6]{DSBL-Bri70} and others.
\begin{DSBL-thm}[Monodromy Theorem]\label{DSBL-1}
The eigenvalues of the monodromy are roots of unity and its Jordan blocks have size at most $(n+1)\times(n+1)$ and size at most $n\times n$ for eigenvalue $1$.
\end{DSBL-thm}

\section{The Gauss-Manin Connection}

Let $\mathrm{M}=\mathrm{M}_s\mathrm{M}_u$ be the decomposition of $\mathrm{M}$ into semisimple part $\mathrm{M}_s$ and unipotent part $\mathrm{M}_u$ and let $\mathrm{N}=-\frac{\log\mathrm{M}_u}{2\pi\mathrm{i}}$ be the nilpotent part of $\mathrm{M}$.
Note that $-2\pi\mathrm{i}\mathrm{N}\in\mathrm{End}_\mathds{Q}(H_\mathds{Q})$ where $H_\mathds{Q}=H\otimes_\mathds{Z}\mathds{Q}$.
Let $H_\mathds{C}=\bigoplus_\lambda H_\mathds{C}^\lambda$ be the decomposition of $H_\mathds{C}=H\otimes_\mathds{Z}\mathds{C}$ into generalized $\lambda$-eigenspaces $H_\mathds{C}^\lambda$ of $\mathrm{M}$ and $\mathrm{M}^\lambda=\mathrm{M}\vert_{H_\mathds{C}^\lambda}$.
Note that $H_\mathds{Q}=H_\mathds{Q}^1\oplus H_\mathds{Q}^{\ne1}$ where $H_\mathds{Q}^1\otimes_\mathds{Q}\mathds{C}=H_\mathds{C}^1$ and $H_\mathds{Q}^{\ne1}\otimes_\mathds{Q}\mathds{C}=\bigoplus_{\lambda\ne1}H_\mathds{C}^\lambda$.
Then there is an inclusion
\[
H_\mathds{C}^{\mathrm{e}^{-2\pi\mathrm{i}\alpha}}\overset{\psi_\alpha}{\longrightarrow}(i_*\mathsc{H})_0
\]
defined by $\psi_\alpha(A)=t^{\alpha+\mathrm{N}}A=t^\alpha\exp(N\log(t))$ with image $C^\alpha=\mathrm{im}(\psi_\alpha)$.
In particular, the operators $\mathrm{M}$ and $\mathrm{N}$ act on $C^\alpha$.
The following lemma is an immediate consequence of the definition of $\psi_\alpha$.
\begin{DSBL-lem}\label{DSBL-6}\
\begin{enumerate}
\item $t\circ\psi_\alpha=\psi_{\alpha+1}$ and $\partial_t\circ\psi_\alpha=\psi_{\alpha-1}\circ(\alpha+\mathrm{N})$.
\item $t:C^\alpha\longrightarrow C^{\alpha+1}$ is bijective and $\partial_t:C^\alpha\longrightarrow C^{\alpha-1}$ is bijective if $\alpha\ne0$.
\item On $C^\alpha$, $t\partial_t-\alpha=\mathrm{N}$ and $\exp(-2\pi\mathrm{i} t\partial_t)=\mathrm{M}^{\mathrm{e}^{-2\pi\mathrm{i}\alpha}}$.
\item $C^\alpha=\ker(t\partial_t-\alpha)^{n+1}$.
\end{enumerate}
\end{DSBL-lem}
\begin{DSBL-dfn}
We call $G=\bigoplus_{-1<\alpha\le0}\mathds{C}\{t\}[t^{-1}] C^\alpha\subset(i_*\mathsc{H})_0$ the local Gauss-Manin connection.
\end{DSBL-dfn}
The local Gauss-Manin connection is a $\mu$-dimensional $\mathds{C}\{t\}[t^{-1}]$-vectorspace and a regular $\mathds{C}\{t\}[\partial_t]$-module.
The generalized $\alpha$-eigenspaces $C^\alpha$ of the operator $t\partial_t$ define the decreasing filtration on $G$ by free $\mathds{C}\{t\}$-modules
\[
V^\alpha=\bigoplus_{\alpha\le\beta<\alpha+1}\mathds{C}\{t\} C^\beta,\quad
V^{>\alpha}=\bigoplus_{\alpha<\beta\le\alpha+1}\mathds{C}\{t\} C^\beta
\]
of rank $\mu$, the V-filtration.
In contrast to the $\psi_\alpha$ and $C^\alpha$, the $V^\alpha$ are independent of the coordinate $t$.
The $C^\alpha$ define a splitting
\[
C^\alpha\cong V^\alpha/V^{>\alpha}=\mathrm{gr}^\alpha_VG
\]
of the V-filtration and we denote by $\mathrm{lead}_V$ the leading term with respect to this splitting.
The ring $\mathds{C}\{t\}$ is a free module of rank $1$ over the ring 
\[
\mathds{C}\{\!\{s\}\!\}=\Bigl\{\sum_{k=0}^\infty a_ks^k\in\mathds{C}[\![s]\!]\Big\vert\sum_{k=0}^\infty\frac{a_k}{k!}t^k\in\mathds{C}\{t\}\Bigr\}
\]
where $s=\int_0^1\mathrm{d} t$ acts by integration.
This fact is generalized by the following lemma \cite[1.3.11]{DSBL-Sch02}.
\begin{DSBL-lem}
The action of $s=\partial_t^{-1}$ on $V^{>-1}$ extends to a $\mathds{C}\{\!\{s\}\!\}$-module structure and $V^{>-1}$ is a free $\mathds{C}\{\!\{s\}\!\}$-module of rank $\mu$.
\end{DSBL-lem}
Since $[\partial_t,t]=1$, $[t,s]=s^2$ and hence
\[
t=s^2\partial_s,\quad\partial_t t=s\partial_s.
\]
We call a free $\mathds{C}\{\!\{s\}\!\}$-submodule of $V^{>-1}$ of rank $\mu$ a $\mathds{C}\{\!\{s\}\!\}$-lattice and call a $t\partial_t$-invariant $\mathds{C}\{\!\{s\}\!\}$-lattice saturated.
A basis $\underline{e}$ of a $\mathds{C}\{\!\{s\}\!\}$-lattice defines a matrix $A=\sum_{k\ge0}A_ks^k$ of $t$ by $t\underline{e}=\underline{e}A$ such that
\[
t\cong A+s^2\partial_s
\]
is the basis representation of $t$.

\section{The Brieskorn Lattice}\label{DSBL-13}

The description of cohomology in terms of holomorphic differential forms by the de Rham isomorphism leads to the definition of the Brieskorn lattice
\[
H''=\Omega_{X,0}^{n+1}/\mathrm{d} f\wedge\mathrm{d}\Omega_{X,0}^{n-1}.
\]
By E.~Brieskorn \cite[1.5]{DSBL-Bri70} and M.~Sebastiani \cite{DSBL-Seb70}, the Brieskorn lattice is the stalk at $0$ of a locally free $\mathsc{O}_T$-sheaf $\mathsc{H}''$ of rank $\mu$ with $\mathsc{H}''\vert_{T'}\cong\mathsc{H}$ and hence $H''\subset(i_*\mathsc{H})_0$.
The regularity of the Gauss-Manin connection proved by E.~Brieskorn \cite[2.2]{DSBL-Bri70} implies that $H''\subset G$.
B.~Malgrange \cite[4.5]{DSBL-Mal74} improved this result by the following theorem.
\begin{DSBL-thm}\label{DSBL-3}
$H''\subset V^{-1}$.
\end{DSBL-thm}
By E.~Brieskorn \cite[1.5]{DSBL-Bri70}, the Leray residue formula can be used to express the action of $\partial_t$ in terms of differential forms by $\partial_t[\mathrm{d} f\wedge\omega]=[\mathrm{d}\omega]$.
In particular, $sH''\subset H''$ and 
\[
H''/sH''\cong\Omega_{X,0}^{n+1}/\mathrm{d} f\wedge\Omega_{X,0}^n\cong\mathds{C}\{\underline{x}\}/\langle\underline{\partial}(f)\rangle.
\]
Since the $V^{>-1}$ is a $\mathds{C}\{\!\{s\}\!\}$-module, theorem \ref{DSBL-3} implies that $H''$ is a free $\mathds{C}\{\!\{s\}\!\}$-module of rank $\mu$ and the action of $s$ can be expressed in terms of differential forms by
\[
s[\mathrm{d}\omega]=[\mathrm{d} f\wedge\omega].
\]
For computational purposes, we may restrict our attention to the completion of the Brieskorn lattice.
E.~Brieskorn \cite[3.4]{DSBL-Bri70} proved the following theorem.
\begin{DSBL-thm}\label{DSBL-7}
The $\mathfrak{m}_{X,0}$- and $\mathfrak{m}_{T,0}$-adic topologies on $H''$ coincide.
\end{DSBL-thm}
While the proof of theorem \ref{DSBL-7} is highly non-trivial, the analogous statement for the $\mathds{C}\{\!\{s\}\!\}$-structure of the Brieskorn lattice is quite elementary \cite[1.5.4]{DSBL-Sch02}.
\begin{DSBL-prp}\label{DSBL-8}
The $\mathfrak{m}_{X,0}$- and $\mathfrak{m}_{\mathds{C}\{\!\{s\}\!\}}$-adic topologies on $H''$ coincide.
\end{DSBL-prp}
We call the completion $\widehat H''$ of $H''$ the formal Brieskorn lattice.
Since completion is faithfully flat, $\widehat H''$ is a free $\mathds{C}[\![s]\!]$-module of rank $\mu$ with a differential operator $t=s^2\partial_s$.
The equality $[\underline{\partial}(f)\overline{g}\mathrm{d}\underline{x}]=s[\underline{\partial}(\overline{g})\mathrm{d}\underline{x}]$ motivates to consider the differential relation $\underline{\partial}(f)-s\underline{\partial}$.
It is not difficult to prove that it defines the formal Brieskorn lattice as a quotient of $\mathds{C}[\![s,\underline{x}]\!]$ \cite[1.5.6]{DSBL-Sch02}.
\begin{DSBL-prp}\label{DSBL-2}
\[
\mathds{C}[\![s,\underline{x}]\!]\overset{\pi_H}{\longrightarrow}\mathds{C}[\![s,\underline{x}]\!]/\langle\underline{\partial}(f)-s\underline{\partial}\rangle\mathds{C}[\![s,\underline{x}]\!]\cong_{\mathds{C}[\![s]\!]}\widehat H''.
\]
\end{DSBL-prp}
Proposition \ref{DSBL-2} is the starting point for an algorithmic approach to the local Gauss-Manin connection.
Let $<_{\underline x}$ be a local degree ordering on $\mathds{C}[\![\underline{x}]\!]$ such that $\deg(\underline{x})<\underline{0}$ and $\deg(\underline{\partial})=-\deg(\underline{x})>\underline{0}$.
One can compute a polynomial standard basis $\underline{g}$ of the Jacobian ideal $\langle\underline{\partial}(f)\rangle$ and a polynomial transformation matrix $B=\bigl(\overline{b}^j\bigr)^j$ such that $\underline{g}=\underline{\partial}(f)B$.
By Nakayama's lemma, $\underline{m}=(\underline{x}^{\underline{\beta}})_{\underline{x}^{\underline{\beta}}\notin\langle\mathrm{lead}(\underline{g})\rangle}$ represents a $\mathds{C}[\![s]\!]$-basis $[\underline{m}]$ of $\widehat H''$.
Let $<_s$ be the local degree ordering on $\mathds{C}[\![s]\!]$ and let $<=(<_s,<_{\underline{x}})$ be the block ordering of $<_s$ and $<_{\underline{x}}$ on $\mathds{C}[\![s,\underline{x}]\!]$.
\begin{DSBL-dfn}\
\begin{enumerate}
\item $\underline{h}=\bigl((g_j-s\underline{\partial}\overline{b}^j)\underline{x}^{\underline \beta}\bigr)_{j,\underline{\beta}}$.
\item $\deg(s)=\min\deg(\underline{m})+2\min\deg(\underline{x})<0$.
\item $N=(N_K)_{K\ge0}$ with $N_K=K\deg(s)-2\min\deg(\underline{x})$.
\item $V=(V_K)_{K\ge0}$ with $V_K=\bigl\{p\in\mathds{C}[\![s,\underline{x}]\!]\big\vert\deg(p)<N_K\bigr\}+\langle s\rangle^K\subset\mathds{C}[\![s,\underline{x}]\!]$.
\end{enumerate}
\end{DSBL-dfn}
Since $\widehat H''$ is a free $\mathds{C}[\![s]\!]$-module, $\underline{h}$ is a standard basis of the $\mathds{C}[\![s]\!]$-module $\langle\underline{\partial}(f)-s\underline{\partial}\rangle\mathds{C}[\![s,\underline{x}]\!]$.
The following lemma is technical but not very deep and can be generalized to formal differential deformations \cite[2.2.10]{DSBL-Sch02}.
\begin{DSBL-lem}\label{DSBL-5}
$V=(V_K)_{K\ge 0}$ is a basis of the $\langle s,\underline{x}\rangle$-adic topology of $\mathds{C}[\![s,\underline{x}]\!]$ with $\pi_H(V_K)=\langle s\rangle^K\widehat H''$.
If $s^\alpha\mathrm{lead}(h_{j,\underline{\beta}})\in V_K$ then $s^\alpha h_{j,\underline{\beta}}\in V_K$.
\end{DSBL-lem}
Lemma \ref{DSBL-5} leads to a normal form algorithm for the Brieskorn lattice \cite[2.2.12]{DSBL-Sch02}.
It computes a normal form with respect to $\underline{h}$ and hence the $[\underline{m}]$-basis representation in $H''$.
The normal form computation up to a given degree can be continued up to any higher degree without additional computational effort.
The normal form algorithm for the Brieskorn lattice is a special case of a modification of Buchberger's normal form algorithm \cite{DSBL-Buc85} for power series rings where termination is replaced by adic convergence \cite[2.1.19]{DSBL-Sch02}.

\section{Mixed Hodge Structure}

By lemma \ref{DSBL-6}, there is a $\mathds{C}$-isomorphism
\[
H_\mathds{C}=\bigoplus_{-1<\alpha\le0}H_\mathds{C}^{\mathrm{e}^{-2\pi\mathrm{i}\alpha}}\overset{\psi}{\longrightarrow}\bigoplus_{-1<\alpha\le0}C^\alpha\cong V^{>-1}/sV^{>-1}
\]
defined by $\psi=\bigoplus_{-1<\alpha\le0}\psi_\alpha$ and the monodromy $\mathrm{M}$ on $H_\mathds{C}$ corresponds to $\exp(-2\pi\mathrm{i} t\partial_t)$ on $\bigoplus_{-1<\alpha\le0}C^\alpha$.

The Hodge filtration $F=(F_k)_{k\in\mathds{Z}}$ on $V^{>-1}$ defined by J.~Scherk and J.H.M.~Steenbrink \cite{DSBL-SS85} is the increasing filtration by the free $\mathds{C}\{\!\{s\}\!\}$-modules
\[
F_k=F^{n-k}=(s^{-k}H'')\cap V^{>-1}
\]
of rank $\mu$.
Via the splitting $C^\alpha\cong\mathrm{gr}_V^\alpha V^{>-1}$, the Hodge filtration induces an increasing Hodge filtration $FC^\alpha$ by $\mathds{C}$-vectorspaces on $C^\alpha$ and, via $\psi$, on $H_\mathds{C}$.
The  nilpotent operator $-2\pi\mathrm{i}\mathrm{N}\in\mathrm{End}_\mathds{Q}(H_\mathds{Q})$ defines an increasing weight filtration $W=(W_k)_{k\in\mathds{Z}}$ centered at $n$ resp. $n+1$ on $H_\mathds{Q}^{\ne1}$ resp. $H_\mathds{Q}^1$.
\begin{DSBL-thm}\label{DSBL-11}
The weight filtration $W$ on $H_\mathds{Q}$ and the Hodge filtration $F$ on $H_\mathds{C}$ define a mixed Hodge structure on the cohomology $H$ of the Milnor fibre and the operator $\mathrm{N}$ is a morphism of mixed Hodge structures of type $(-1,-1)$.
\end{DSBL-thm}
The mixed Hodge structure on the cohomology of the Milnor fibre was discovered by J.H.M.~Steenbrink \cite{DSBL-Ste76} and described in terms of the Brieskorn lattice by A.N.~Varchenko \cite{DSBL-Var82a}.

The nilpotent operator $N$ on $C^\alpha$ defines an increasing weight filtration $W=(W_k)_{k\in\mathds{Z}}$ centered at $n$ on $C^\alpha$.
By definition $\mathrm{N}$ commutes with $\psi_\alpha$ and hence
\[
\psi_\alpha\bigl(WH_\mathds{C}^{\mathrm{e}^{-2\pi\mathrm{i}\alpha}}\bigr)=
\begin{cases}
WC^\alpha, & \alpha\notin\mathds{Z},\\
W[-1]C^\alpha, & \alpha\in\mathds{Z}.
\end{cases}
\]
The weight filtration $W=\bigoplus_{-1<\alpha\le0}\mathds{C}\{\!\{s\}\!\} WC^\alpha$ on $V^{>-1}$ by free $\mathds{C}\{\!\{s\}\!\}$-modules induces $WC^\alpha$ via the splitting $C^\alpha\cong\mathrm{gr}_V^\alpha V^{>-1}$. 

The spectral pairs are those pairs $(\alpha,l)\in\mathds{Q}\times\mathds{Z}$ with positive multiplicity
\[
d^\alpha_l=\dim_\mathds{C}\mathrm{gr}^W_l\mathrm{gr}_V^\alpha\mathrm{gr}^F_0V^{>-1}.
\]
Via the isomorphism $\psi$, they correspond to the Hodge numbers
\[
h^{p,l-p}_\lambda=\dim_\mathds{C}\mathrm{gr}_F^p\mathrm{gr}^W_l H_\mathds{C}^\lambda
\]
by $d^{\alpha+p}_l=h_{\mathrm{e}^{-2\pi\mathrm{i}\alpha}}^{n-p,l-n+p}$ for $-1<\alpha<0$ and $d^p_l=h_1^{n-p,l+1-n+p}$ and inherit the symmetry properties
\[
d^\alpha_l=d^{2n-l-1-\alpha}_l,\quad d^\alpha_l=d^{\alpha-n+l}_{2n-l},\quad d^\alpha_l=d^{n-1-\alpha}_{2n-l}
\]
from the mixed Hodge structure.
The spectral numbers are those numbers $\alpha\in\mathds{Q}$ with positive multiplicity
\[
d^\alpha=\dim_\mathds{C}\mathrm{gr}_V^\alpha\mathrm{gr}^F_0V^{>-1}=\sum_{l\in\mathds{Z}}d^\alpha_l
\]
and have the symmetry property $d^\alpha=d^{n-1-\alpha}$.

\section{M. Saito's Basis}

By P.~Deligne \cite[1.2.8]{DSBL-Del72}, a morphism of mixed Hodge structures is strict for the Hodge filtration.
In particular, by theorem \ref{DSBL-11}, $\mathrm{N}$ is strict for the Hodge filtration on $H_\mathds{C}$ and on $\mathrm{gr}_V V^{>-1}$.
Hence, there is a direct sum decomposition $F_kC^\alpha=\bigoplus_{j\le k}C^{\alpha,j}$ such that $\mathrm{N}(C^{\alpha,k})\subset C^{\alpha,k+1}$, and $sC^{\alpha,k}\subset C^{\alpha+1,k-1}$.
By definition of the Hodge filtration,
\[
\mathrm{lead}_V(H'')=\sum_{\alpha\in\mathds{Q}}\sum_{k\le0}\mathds{C}\{\!\{s\}\!\} C^{\alpha,k}=\bigoplus_{\alpha\in\mathds{Q}}\mathds{C}\{\!\{s\}\!\} G^\alpha
\]
where $G^\alpha=C^{\alpha,0}$.
Let $<_{\mathds{Q}\times\mathds{Z}}=(>_\mathds{Q},>_\mathds{Z})$ be the block ordering of $>_\mathds{Q}$ and $>_\mathds{Z}$ on the index set $\mathds{Q}\times\mathds{Z}$.
Then the Hodge filtration defines a refinement of the V-filtration on $V^{>-1}$ by free $\mathds{C}\{\!\{s\}\!\}$-modules $V^{\alpha,k}=F_kC^\alpha\oplus V^{>\alpha}$ of rank $\mu$ and the $C^{\alpha,k}$ define a splitting of this refined filtration compatible with $s$.
We call the refinement the Hodge refinement and the splitting a Hodge splitting.
The following lemma follows essentially from the fact that $\mathds{C}\{\!\{s\}\!\}$ is a discrete valuation ring \cite[1.10.5,1.10.10]{DSBL-Sch02}.
\begin{DSBL-lem}
Let $H$ be a $\mathds{C}\{\!\{s\}\!\}$-lattice and $C^{\alpha,k}$ a splitting of a refinement of the V-filtration compatible with $s$.
Then a minimal standard basis of $H$ is a $\mathds{C}\{\!\{s\}\!\}$-basis and there is a reduced minimal standard basis of $H$.
\end{DSBL-lem}
In particular, there is a reduced minimal standard basis of $H''$ for a Hodge splitting.
The following proposition follows essentially from lemma \ref{DSBL-6}.3 \cite[1.10.12]{DSBL-Sch02}.
\begin{DSBL-prp}\label{DSBL-4}
Let $\underline{h}$ be a reduced minimal standard basis of $H''$ for a Hodge splitting.
Then the $\underline{h}$-matrix $A$ of $t$ has degree $1$.
In particular, 
\[
\bigl(H'',t\bigr)\overset{\underline{h}}{\longleftarrow}\bigl(\mathds{C}\{\!\{s\}\!\}^\mu,A_0+A_1s+s^2\partial_s\bigr)
\]
is an isomorphism.
Moreover, $A_1$ is semisimple with eigenvalues the spectral numbers of $f$ added by $1$ and $\mathrm{gr}_V(A_0)$ can be identified with $\mathrm{N}$.
\end{DSBL-prp}
Note that the matrices $A_0$ and $A_1$ in proposition \ref{DSBL-4} determine the differential structure of the Brieskorn lattice.
M. Saito \cite{DSBL-Sai89} first constructed a $\mathds{C}\{\!\{s\}\!\}$-basis of $H''$ as in proposition \ref{DSBL-4} without calling it a reduced minimal standard basis.

\section{The Algorithm}\label{DSBL-12}

We describe an algorithm to compute $A_0$ and $A_1$ as in proposition \ref{DSBL-4} \cite{DSBL-Sch02}.
This algorithm can be simplified to compute the complex monodromy, the spectral numbers, or the spectral pairs only \cite{DSBL-Sch02}.

The normal form algorithm for the Brieskorn lattice in section \ref{DSBL-13} computes the $[\underline{m}]$-matrix $A=\sum_{k\ge0}A_ks^k$ of $t$ defined by $t[\underline{m}]=[f\underline{m}]=[\underline{m}]A$ up to any degree.
We identify the columns of a matrix $H$ with the generators of a submodule $\langle H\rangle\subset\mathds{C}[\![s]\!]^\mu$ and denote by $E$ the unit matrix.
Then $\langle E\rangle$ is the $[\underline{m}]$-basis representation of $\widehat H''$.
Hence, the following two statements hold for $\underline{h}=[\underline{m}]$ with $\kappa=0$ and $H=E$.
\begin{enumerate}
\item[($H_{\underline{h}}$)] One can compute $\kappa\ge 0$ and a $\mu\times\mu$-matrix $H$ with coefficients in $\mathds{C}[s]$ of degree at most $\kappa$ such that $\langle H\rangle$ is the $\underline{h}$-basis representation of $\widehat H''$ and $s^\kappa\langle E\rangle\subset\langle H\rangle$.
\item[($A_{\underline{h}}$)] One can compute the $\underline{h}$-matrix $A$ of $t$ up to any degree.
\end{enumerate}
Step by step, we improve the $\mathds{C}[\![s]\!]$-basis $\underline{h}$ and show that $(H_{\underline{h}})$ and $(A_{\underline{h}})$ hold.
After the last step, $A_0$ and $A_1$ as in proposition \ref{DSBL-4} can be computed by a basis transformation of $A$ to a reduced minimal standard basis of $\langle H\rangle$ up to a certain degree bound.

We call the canonical projection $\mathrm{jet}_k:\mathds{C}[\![s]\!]\longrightarrow \bigoplus_{j=0}^k\mathds{C} s^j$ the $k$-jet.
Let the monomial ordering on $\mathds{C}[\![s]\!]^\mu=\mathds{C}[\![s]\!]\otimes_\mathds{C}\mathds{C}^\mu$ be the block ordering $<=(<_s,>_\mu)$ of the local degree ordering $<_s$ on $\mathds{C}[\![s]\!]$ and the inverse ordering $>_\mu$ on the indices of the basis elements of $\mathds{C}^\mu$.

\subsection{The Saturation of $H''$}
In this step, we show that $(H_{\underline{h}})$ and $(A_{\underline{h}})$ hold for a $\mathds{C}[\![s]\!]$-basis $\underline{h}$ of a saturated $\mathds{C}[\![s]\!]$-lattice.

The increasing sequence of $\mathds{C}[\![s]\!]$-lattices defined by
\[
\widehat H''_0=\widehat H'',\quad
\widehat H''_{k+1}=s\widehat H''_{k}+t\widehat H''_{k}\subset \widehat H''
\]
is stationary since $\widehat H''$ is noetherian.
Hence, the saturation $\widehat H''_\infty=\bigcup_{k\ge0}\widehat H''_k$ of $\widehat H''$ is a saturated $\mathds{C}[\![s]\!]$-lattice.
The $[\underline{m}]$-basis representation $\langle H_k\rangle$ of $\widehat H''_k$ can be computed by
\[
H_0=Q_{-1}=E,\quad
Q_k=\bigl(\mathrm{jet}_k(A)+s^2\partial_s\bigr)Q_{k-1},\quad
H_{k+1}=(sH_k\vert Q_k).
\]
We successively compute the $H_k$ and check in each step if $\langle Q_k\rangle\subset\langle H_k\rangle$ by a standard basis and normal form computation.
If $\langle Q_k\rangle\subset\langle H_k\rangle$ then we stop the computation and set $\kappa=k$ and $H_\infty=H_\kappa$.
Then $\langle H_\infty\rangle$ is the $[\underline{m}]$-basis representation of $\widehat H''_\infty$.
We replace $H_\infty$ by a minimal standard basis of $\langle H_\infty\rangle$.
Then $\underline{h}=s^{-\kappa}\underline{h}H_\infty$ is a $\mathds{C}[\![s]\!]$-basis of a saturated $\mathds{C}[\![s]\!]$-lattice.
By a normal form computation with respect to $H_\infty$ up to degree $\kappa$, we compute the $\underline{h}$-basis representation $\langle H_\infty^{-1} s^\kappa E\rangle=\langle\mathrm{jet}_\kappa(H_\infty^{-1} s^\kappa E)\rangle$ of $\widehat H''$.
Since $\langle H_\infty\rangle\subset\langle E\rangle$, $s^\kappa\langle E\rangle\subset\langle H_\infty^{-1} s^\kappa E\rangle$.
By a normal form computation with respect to $H_\infty$ up to degree $\kappa+k$, one can compute the $k$-jet 
\[
\mathrm{jet}_k\bigl(H_\infty^{-1}(A-\kappa sE+s^2\partial_s)H_\infty\bigr)=\mathrm{jet}_k\bigl(H_\infty^{-1}\bigl(\mathrm{jet}_{\kappa+k}(A-\kappa sE)+s^2\partial_s\bigr)H_\infty\bigr)
\]
of the $\underline{h}$-matrix of $t$ for any $k\ge0$.

\subsection{The V-Filtration}
In this step, we show that $(H_{\underline{v}})$ and $(A_{\underline{v}})$ hold for a $<_\mathds{Q}$-increasingly ordered $\mathds{C}[\![s]\!]$-basis $\underline{v}$ of a $\widehat V^\alpha$ compatible with the direct sum decomposition $\widehat V^\alpha/s\widehat V^\alpha\cong\bigoplus_{\alpha\le\beta<\alpha+1}C^\beta$.

Since $\underline{h}$ is a $\mathds{C}[\![s]\!]$-basis of a saturated $\mathds{C}\{\!\{s\}\!\}$-lattice, $A_0=0$ and, by theorem \ref{DSBL-1}, the eigenvalues of $A_1$ are rational.
In order to compute the eigenvalues of $A_1$, we transform $A_1$ to Hessenberg form and factorize the characteristic polynomials of its blocks.
Then we compute a constant $\mathds{C}[\![s]\!]$-basis transformation such that $A_1=\mathrm{diag}(\alpha_1,\dots,\alpha_\mu)+N$ with $\alpha_1\le\cdots\le\alpha_\mu$ where $\mathrm{diag}(\alpha_1,\dots,\alpha_\mu)$ denotes the diagonal matrix with entries $\alpha_1,\dots,\alpha_\mu$.
If $\alpha_\mu-\alpha_1<1$ then $\underline{v}=\underline{h}$ is a $<_\mathds{Q}$-increasingly ordered $\mathds{C}[\![s]\!]$-basis $\underline{v}$ of a $\widehat V^\alpha$ compatible with the direct sum decomposition $\widehat V^\alpha/s\widehat V^\alpha\cong\bigoplus_{\alpha\le\beta<\alpha+1}C^\beta$.
If $\alpha_\mu-\alpha_1\ge1$ then we proceed as follows.
Let
\[
A=\begin{pmatrix}A^{1,1}&A^{1,2}\\A^{2,1}&A^{2,2}\end{pmatrix}\index{Aij@$A^{i,j}$}
\]
such that $A_0=0$, $A^{1,2}_1=0$, $A^{2,1}_1=0$, and the eigenvalues of $A^{1,1}_1$ are the eigenvalues $\alpha$ of $A_1$ with $\alpha<\alpha_1+1$.
Then the $\mathds{C}[\![s]\!][s^{-1}]$-basis transformation
{\small
\[
H\mapsto\begin{pmatrix}\frac{1}{s}&0\\0&1\end{pmatrix}H,\quad
A\mapsto\begin{pmatrix}\frac{1}{s}&0\\0&1\end{pmatrix}\bigl(A+s^2\partial_s\bigr)\begin{pmatrix}s&0\\0&1\end{pmatrix}
=\begin{pmatrix}A^{1,1}+s&\frac{1}{s}A^{1,2}\\sA^{2,1}&A^{2,2}\end{pmatrix}
\]
}
decreases $\alpha_\mu-\alpha_1$ and the degree up to which $A$ is computed by $1$ and increases $\kappa$ by $1$.
After at most $n$ such transformations, $\alpha_\mu-\alpha_1<1$.

\subsection{The Canonical V-Splitting}
In this step, we show that $(H_{\underline{c}})$ and $(A_{\underline{c}})$ hold for a $<_\mathds{Q}$-increasingly ordered $\mathds{C}$-basis $\underline{c}$ of a direct sum $\bigoplus_{\alpha\le\beta<\alpha+1}C^\beta$ compatible with the direct sum.

Let $\underline{c}$ be the image of $[\underline{v}]$ under the splitting $\widehat V^\alpha/s\widehat V^\alpha\cong\bigoplus_{\alpha\le\beta<\alpha+1}C^\beta$.
By Nakayama's lemma, $\underline{c}$ is a $\mathds{C}$-basis of $\bigoplus_{\alpha\le\beta<\alpha+1}C^\beta$ compatible with the direct sum.
The eigenvalues of the commutator $[\cdot,A_1]\in\mathrm{End}_\mathds{C}(\mathds{C}^{\mu^2})$ are the differences of the eigenvalues of $A_1$.
Since $\alpha_\mu-\alpha_1<1$, $[\cdot,A_1]-k\in\mathrm{GL}_{\mu^2}(\mathds{C})$ for $k\ge1$.
Let $U=\sum_{j=0}^\infty U_js^j$ be the $\mathds{C}[\![s]\!]$-basis transformation defined by $\underline{c}=\underline{v}U$.
Then $U_0=E$ and $UA_1s=(A+s^2\partial_s)U$ or equivalently
\[
U_k=\bigl([\cdot,A_1]-k\bigr)^{-1}\sum_{j=0}^{k-1}A_{k-j+1}U_j
\]
for $k\ge1$ and hence one can compute $U$ up to any degree.
Since $U_0=E$ and $\kappa\ge0$, $\mathrm{jet}_\kappa(U)$ is a minimal standard basis of $\langle E\rangle$.
By a normal form computation with respect to $U$ up to degree $\kappa$, we compute the $\underline{c}$-basis representation $\langle U^{-1}H\rangle=\langle\mathrm{jet}_\kappa(\mathrm{jet}_\kappa(U)^{-1}H)\rangle$ of $\widehat H''$ and $A_1s$ is the $\underline{c}$-matrix of $t$.

\subsection{A Hodge Splitting}\ \\
In this step, we show that $(H_{\underline{f}})$ and $(A_{\underline{f}})$ hold for a $<_{\mathds{Q}\times\mathds{Z}}$-decreasingly ordered $\mathds{C}$-basis $\underline{f}$ of a direct sum $\bigoplus_{\alpha\le\beta<\alpha+1}\bigoplus_{k\in\mathds{Z}}C^{\beta,k}$ compatible with the direct sum and that one can compute $A_0$ and $A_1$ as in proposition \ref{DSBL-4}.

We compute a standard basis of $H$ up to degree $\kappa$ in order to compute the $\underline{c}$-basis representation of the Hodge filtration $F$.
The nilpotent part of $A_1$ is the $\underline{c}$-basis representation of $\mathrm{N}$.
By computing images and quotients of $\mathds{C}$-vectorspaces, we compute the $\underline{c}$-basis representation of a Hodge splitting $F_kC^\beta=\bigoplus_{j\le k}C^{\beta,j}$.
Then we compute a constant $\mathds{C}[\![s]\!]$-basis transformation $\underline{f}=\underline{c}U$ such that $\underline{f}$ is a $<_{\mathds{Q}\times\mathds{Z}}$-decreasingly ordered $\mathds{C}$-basis of the direct sum $\bigoplus_{\alpha\le\beta<\alpha+1}\bigoplus_{k\in\mathds{Z}}C^{\beta,k}$ compatible with the direct sum.

We replace $H$ by a reduced minimal standard basis of $\langle H\rangle$ up to degree $\kappa+1$.
By a normal form computation with respect to $H$ up to degree $\kappa+1$, we compute the $1$-jet 
\[
\mathrm{jet}_1\bigl(H^{-1}(A+s^2\partial_s)H\bigr)=\mathrm{jet}_1\bigl(\mathrm{jet}_{\kappa+1}(H)^{-1}\bigl(\mathrm{jet}_{\kappa+1}(A)+s^2\partial_s\bigr)\mathrm{jet}_{\kappa+1}(H)\bigr)
\]
of the $\underline{c}H$-matrix $A$ of $t$ in order to compute $A_0$ and $A_1$ as in proposition \ref{DSBL-5}.

\section{An Example}

The algorithm in section \ref{DSBL-12} is implemented in the computer algebra system {\sc Singular} \cite{DSBL-GPS02} in the procedure {\tt tmatrix} in the library {\tt gaussman.lib} \cite{DSBL-Sch01b}.
In an example {\sc Singular} session, we compute the differential structure of the Brieskorn lattice of the singularity of type $T_{2,5,5}$ defined by the polynomial $f=x^2y^2+x^5+y^5$.

First, we load the {\sc Singular} library {\tt gaussman.lib}:
{\small\begin{verbatim}
> LIB "gaussman.lib";
\end{verbatim}}
Then, we define the local ring $R=\mathds{Q}[x,y]_{\langle x,y\rangle}$ with the local degree ordering {\tt ds} as monomial ordering and the polynomial $f=x^2y^2+x^5+y^5\in R$:
{\small\begin{verbatim}
> ring R=0,(x,y),ds;
> poly f=x2y2+x5+y5;
\end{verbatim}}
Finally, we compute $A_0$ and $A_1$ as in proposition \ref{DSBL-5}:
{\small\begin{verbatim}
> list A=tmatrix(f);
\end{verbatim}}
The result is the list {\tt A=A[1],A[2]} such that ${\tt A[i+1]}=A_i$ and 
{\small
\[
\setcounter{MaxMatrixCols}{20}
A_0=
\begin{pmatrix}
0&0&\cdots&0\\
\vdots&\vdots&&\vdots\\
0&0&\cdots&0\\
1&0&\cdots&0
\end{pmatrix},\quad
A_1=\mathrm{diag}\Bigl(\frac{1}{2},\frac{7}{10},\frac{7}{10},\frac{9}{10},\frac{9}{10},1,\frac{11}{10},\frac{11}{10},\frac{13}{10},\frac{13}{10},\frac{3}{2}\Bigr).
\]
}
By proposition \ref{DSBL-5}, $\bigl(H'',t\bigr)\cong\bigl(\mathds{C}\{\!\{s\}\!\}^\mu,A_0+sA_1+s^2\partial_s\bigr)$ and the spectral pairs are $\bigl(-\frac{1}{2},2\bigr)$, $\bigl(-\frac{3}{10},1\bigr)^2$, $\bigl(-\frac{1}{10},1\bigr)^2$, $(0,1)$, $\bigl(\frac{1}{10},1\bigr)^2$, $\bigl(\frac{3}{10},1\bigr)^2$, $\bigl(\frac{1}{2},0\bigr)$.

\end{document}